\documentclass[12pt]{amsart}
\usepackage[utf8]{inputenc}
\usepackage[T1]{fontenc}
\usepackage[top=3cm, bottom=3cm, left=3cm, right=3cm]{geometry}
\usepackage{csquotes, amsfonts, amsmath, amsthm, amssymb, stmaryrd, dsfont, enumitem, xcolor, url, graphicx, multicol, float, setspace, hyperref, tikz, multirow, esvect, colortbl, pifont, doi, url}
\usepackage[normalem]{ulem}
\usepackage[breakable, skins]{tcolorbox}

\setlist[itemize]{leftmargin=1cm}
\setlist[enumerate]{leftmargin=1cm}

\theoremstyle{plain}
\newtheorem{Theorem}{Theorem}[section]
\newtheorem{Proposition}[Theorem]{Proposition}

\newtheorem{Corollary}[Theorem]{Corollary}

\theoremstyle{definition}
\newtheorem{Definition}[Theorem]{Definition}
\newtheorem{Remark}[Theorem]{Remark}
\newtheorem{Example}[Theorem]{Example}

\hypersetup{
    colorlinks=true,
    linkcolor=black,
    urlcolor=black,
    citecolor=black
}

\newcommand{\K}{\mathbb{K}}
\newcommand{\R}{\mathbb{R}}
\newcommand{\C}{\mathbb{C}}
\newcommand{\D}{\mathbb{D}}
\newcommand{\Z}{\mathbb{Z}}
\newcommand{\T}{\mathbb{T}}

\newcommand{\dd}{\mathrm{d}}

\newcommand{\Hol}{\mathop{\rm Hol}\nolimits}

\newcommand{\indic}{\mathds{1}}

\newcommand{\Cont}{\mathcal{C}}

\newcommand{\norm}[1]{\left\Vert#1\right\Vert}

\newcommand{\abs}[1]{\left\lvert #1 \right\rvert}

\newcommand{\sbt}{\,\begin{picture}(-1,1)(-1,-3)\circle*{3}\end{picture}\ }

\renewcommand{\textbf}[1]{\begingroup\bfseries\mathversion{bold}#1\endgroup}

\title[Distances and seminorms on Fréchet spaces]
{On the relation between distances and seminorms on Fréchet spaces, with application to isometries}

\author{I. Chalendar}
\address{Isabelle CHALENDAR, Université Gustave Eiffel, LAMA, (UMR 8050), 
    UPEM, UPEC, CNRS, F-77454, Marne-la-Vallée (France)}
\email{isabelle.chalendar@univ-eiffel.fr}

 \author{L. Oger}
\address{Lucas OGER, Université Gustave Eiffel, LAMA, (UMR 8050), 
    UPEM, UPEC, CNRS, F-77454, Marne-la-Vallée (France)}
\email{lucas.oger@univ-eiffel.fr}

\author{J. R. Partington}
\address{Jonathan R. PARTINGTON, School of Mathematics, University of Leeds, Leeds LS2 9JT, Yorkshire, U.K.}
\email{j.r.partington@leeds.ac.uk}


\begin{document}

\begin{abstract}
A study is made of linear isometries on Fréchet spaces for which the metric is given in terms of a sequence of seminorms. This establishes sufficient conditions on the growth of the function that defines the metric in terms of the seminorms to ensure that a linear operator preserving the metric also preserves each of these seminorms. As an application, characterizations are given of the isometries on various spaces including those of holomorphic functions on complex domains and continuous functions on open sets, extending the Banach--Stone theorem to surjective and nonsurjective cases.
\end{abstract}

\maketitle

\section{Introduction}

Let $(E, d)$ be a metric vector space. A \emph{linear isometry} is a linear operator $T$ such that
\begin{equation}\label{Eq - Def linear isometry}
    \forall x, y \in E, \quad d(Tx,Ty) = d(x,y).
\end{equation}

Studying the linear isometries of a space allows us to appreciate its geometry, and many applications exist, for instance in linear dynamics \cite{Ansari-Bourdon, Oger}.
In Banach spaces, there is a very convenient way of describing these operators, using the norm $\norm{\cdot}$ endowed with $E$. In this case, using linearity, \eqref{Eq - Def linear isometry} is equivalent to the following property.
\begin{equation}\label{Eq - Def 2 linear isometry}
    \forall x \in E, \quad \norm{Tx} = \norm{x}.
\end{equation}

This research is focused on a generalisation of this property \eqref{Eq - Def 2 linear isometry} in the context of Fréchet spaces, which are not complete for any norm. Instead, we endow $E$ with an increasing family of seminorms, denoted by $\norm{\cdot}_n$ $(n \ge 0)$. This gives the possibility to construct a zoology of distances. Indeed, let $\theta : \R^+ \to \R^+$ be an increasing and subadditive map, such that $\theta(0) = 0$ and $\theta(t) \to 1$ as $t \to +\infty$. Let $r = (r_n)_{n \ge 0}$ be a summable positive sequence. Without loss of generality, we will assume in the following that $\sum r_n = 1$. Set
\[ d_{(\theta,r)}(x,y) := \sum_{n \ge 0} r_n \theta(\norm{x-y}_n),
    \quad x,y \in E. \]
The hypotheses on $\theta$ and $r$ mean that $d_{(\theta,r)}$ is a distance on the space $E$. Moreover, once again by linearity, $T$ is a linear isometry of $E$ if and only if  for all $x \in E$, $d_{(\theta,r)}(Tx,0) = d_{(\theta,r)}(x,0)$.
A natural question is a description of linear isometries for $d_{(\theta,r)}$, only using the seminorms. In other words, if $\mathcal L(E)$ denotes the space of all continuous linear maps $T : E \to E$, which maps $\theta$ and sequences $r$ satisfy the following property?
\begin{equation}\label{Propriété}\tag{$P$}
    \forall T \in \mathcal L(E), \; \forall x \in E, \;
    d_{(\theta,r)}(Tx,0) = d_{(\theta,r)}(x,0) \iff 
    (\forall n \ge 0, \; \norm{Tx}_n = \norm{x}_n).
\end{equation}

It is clear that $(\Longleftarrow)$ is always satisfied. Hence, the main problem is to obtain the reverse implication $(\Longrightarrow)$. We begin by giving a sufficient condition on $\theta$ and $r$ to obtain the property \eqref{Propriété}. Then, we focus on the description of linear isometries on the space $\Hol(\D)$ of all holomorphic functions on the unit disc $\D$, and the space $\Cont(U)$ of all continuous maps on an open set $U$. Hence, we link the study with the theory of weighted composition operators, and Banach-Stone theorem.

%%%%%%%%%%%%%%%%%%%%%%%%%%%%%%%%%%%%%%%%%%
\section{A sufficient condition to satisfy \eqref{Propriété}}
\label{Section 1}

Before stating the main theorem, let us recall the definition of an absolutely continuous map. In the following, we denote $a(t) \lesssim b(t)$ when $a(t) \le C b(t)$ for all $t$, with $C > 0$ a constant independent from $t$.

\begin{Definition}
    A map $\theta : I \subset \R \to \R$ is \textbf{absolutely continuous} if $\theta$ has a derivative $\theta'$ almost everywhere which is integrable in $I$.
\end{Definition}

\begin{Theorem}\label{Thm - Propriété (P)}
    Let $\theta : \R^+ \to \R^+$ be an increasing and subadditive map, such that $\theta(0) = 0$ and $\theta(x) \to 1$ as $x \to +\infty$. Assume that $\theta$ is \textbf{absolutely continuous}, and there exists $\alpha, m, M > 0$ such that
    \[ t < m \implies \theta(t) \lesssim t^\alpha, \qquad
        t > M \implies 1 - \theta(t) \lesssim t^{-\alpha}. \]
    Then, for all summable and positive sequence $r = (r_n)_{n \ge 0}$, and for all increasing positive sequences $a = (a_n)_{n \ge 0}$, $b = (b_n)_{n \ge 0}$, we have
    \[ \left(\forall t > 0, \quad \sum_{n \ge 0} r_n \theta(t a_n) 
        = \sum_{n \ge 0} r_n \theta(t b_n) \right)
        \implies a = b. \]
\end{Theorem}

\begin{proof}
    First, note that since $\theta$ is absolutely continuous, the map $\theta'$ is integrable on $\R^+$. Hence, for all $\rho > 1$, we can define the Frullani integral
    \begin{align*}
        \int_0^{+\infty} \frac{\theta(\rho x) - \theta(x)}{x} \; \dd x
        = \int_0^{+\infty} \frac 1x \left( \int_{x}^{\rho x} \theta'(t) \; \dd t \right) \; \dd x
        & = \int_0^{+\infty} \theta'(t) \left( \int_{t/\rho}^{t} \frac 1x \; \dd x \right) \; \dd t \\
        & = \int_0^{+\infty} \ln(\rho) \theta'(t) \; \dd t < +\infty.
    \end{align*}

    Moreover, if we denote by $\delta_x$ the Dirac measure at $x$, then for all $t > 0$,
    \begin{equation}\label{Eq Et}\tag{$E_t$}
        \sum_{n \ge 0} r_n \theta(t a_n) = \sum_{n \ge 0} r_n \theta(t b_n)
        \iff \int_0^{+\infty} \theta(tx) \; \dd \mu_a(x)
            = \int_0^{+\infty} \theta(tx) \; \dd \mu_b(x),
    \end{equation} 
    with $\mu_a$ and $\mu_b$ the Borel measures defined by
    \[ \mu_a = \sum_{n \ge 0} r_n \delta_{a_n}, \quad
        \mu_b = \sum_{n \ge 0} r_n \delta_{b_n}. \]
    We need to show that if $(E_t)$ is satisfied for all $t > 0$, then $\mu_a = \mu_b$. \medskip

    Set $\rho = e^u$, $x = e^y$ and $F(w) = \theta(e^w)$, for $w \in \R$. Using the fact that the exponential function is a continuous bijection from $\R$ to $\R^+$, and defining the measure $\nu$, for each Borel subset $E$ of $\R$, by $\nu(E)= \mu(\exp(E))$, we get
    \[ \int_0^{+\infty} \theta(\rho x) \; \dd \mu(x)
        = \int_{-\infty}^{+\infty} \theta(e^u e^y) \; \dd \nu(y)
        = \int_{-\infty}^{+\infty} F(u+y) \; \dd \nu(y). \]
        
    We fix $u > 0$. Denoting $G_u(y) = F(u+y) - F(y)$, we obtain $G_u \ge 0$ (since $F$ is increasing), and by the first step of this proof,
    \[ \int_{-\infty}^{+\infty} G_u(y) \; \dd y
        = \int_{-\infty}^{+\infty} [F(u+y)-F(y)] \; \dd y
        = \int_0^{+\infty} \frac{\theta(\rho x) - \theta(x)}{x} \; \dd x < +\infty, \]
    Therefore, $G_u \in L^1(\R)$. \medskip

    Let us consider the map
    \[ H_u(s) = \int_{-\infty}^{+\infty} G_u(y+s) \; \dd \nu(y). \]
    Then, its Fourier transform is
    \begin{align*}
        \widehat{H_u}(z)
        & = \int_{-\infty}^{+\infty} \int_{-\infty}^{+\infty}
            G_u(y+s) e^{-isz} \; \dd s \; \dd \nu(y) \\
        & = \int_{-\infty}^{+\infty} \int_{-\infty}^{+\infty}
            G_u(w) e^{-i(w-y)z} \; \dd w \; \dd \nu(y) \\
        & = \left(\int_{-\infty}^{+\infty} G_u(w) e^{-iwz} \; \dd w\right)
            \left(\int_{-\infty}^{+\infty} e^{iyz} \; \dd \nu(y)\right)
        = \widehat{G_u}(z) \times \widehat{\nu}(-z).
    \end{align*}
    We can therefore uniquely determine $\widehat{\nu}(z)$, provided that $\widehat{G_u}(-z) \neq 0$ for a $u > 0$. If there exists $\alpha > 0$ such that $\theta(t) = O(t^\alpha)$ when $t \to 0$ and $1 - \theta(t) = O(t^{-\alpha})$ when $t \to +\infty$, then $G_u(w) = O(e^{-\alpha w})$ in the neighbourhood of $-\infty$, and $G_u(w) = O(e^{-\alpha w})$ in the neighbourhood of $+\infty$. Thus, $\widehat{G_u}(z)$ is well defined and holomorphic on $\{z \in \C : \abs{\Im(z)} < \alpha\}$. We deduce that the zeros of $\widehat{G_u}(z)$ are isolated, so there exists a dense set where $\widehat \nu (z)$ is uniquely determined. By injectivity of the Fourier transform, the same applies to $\nu$. \medskip

    Finally, for all $t > 0$,
    \[ \int_0^{+\infty} \theta(tx) \; \dd \mu_a(x) - \int_0^{+\infty} \theta(tx) \; \dd \mu_b(x)
        = \int_0^{+\infty} \theta(tx) \; \dd (\mu_a-\mu_b)(x) = 0, \]
    so that $\mu_a - \mu_b = 0$ (since the measure is unique). Thus, $a = b$.
\end{proof}

\begin{Corollary}\label{Coro - Propriété (P)}
   Let $\theta : \R^+ \to \R^+$ be an increasing and subadditive map, such that $\theta(0) = 0$ and $\theta(x) \to 1$ as $x \to +\infty$. Assume that $\theta$ is \textbf{absolutely continuous}, and there exists $\alpha, m, M > 0$ such that
    \[ t < m \implies \theta(t) \lesssim t^\alpha, \qquad
        t > M \implies 1 - \theta(t) \lesssim t^{-\alpha}. \]
    Let $r = (r_n)_{n \ge 0}$ be a summable positive sequence. Then, $\theta$ and $r$ satisfy $(P)$.
\end{Corollary}

\begin{proof}
    Let $T \in \mathcal L(E)$, and $x \in E$ such that $d_{(\theta,r)}(Tx, 0) = d_{(\theta,r)}(x, 0)$.
    For all $n \ge 0$, let us denote $a_n = \norm{Tx}_n$ and $b_n = \norm{x}_n$, and assume that $a_n > 0$ for $n \ge n_1$ and $b_n > 0$ for $n \ge n_2$. Since the seminorms are homogeneous, for all $t > 0$,
    \[ \sum_{n \ge n_1} r_n \theta(ta_n) 
        = \sum_{n \ge n_2} r_n \theta(tb_n). \]
    Then, as $t \to +\infty$, we obtain
    \[ \sum_{n \ge n_1} r_n = \sum_{n \ge n_2} r_n. \]
    Finally, $n_1 = n_2$ because all the $r_n$ are positive. Without loss of generality, assume that $n_1 = n_2 = 0$. By Theorem \ref{Thm - Propriété (P)}, we conclude that for all $n \ge 0$, $a_n = b_n$, that is
    \[ \norm{Tx}_n = \norm{x}_n. \qedhere \]
\end{proof}

\begin{figure}[H]
    \centering
    \includegraphics[width=0.5\linewidth]{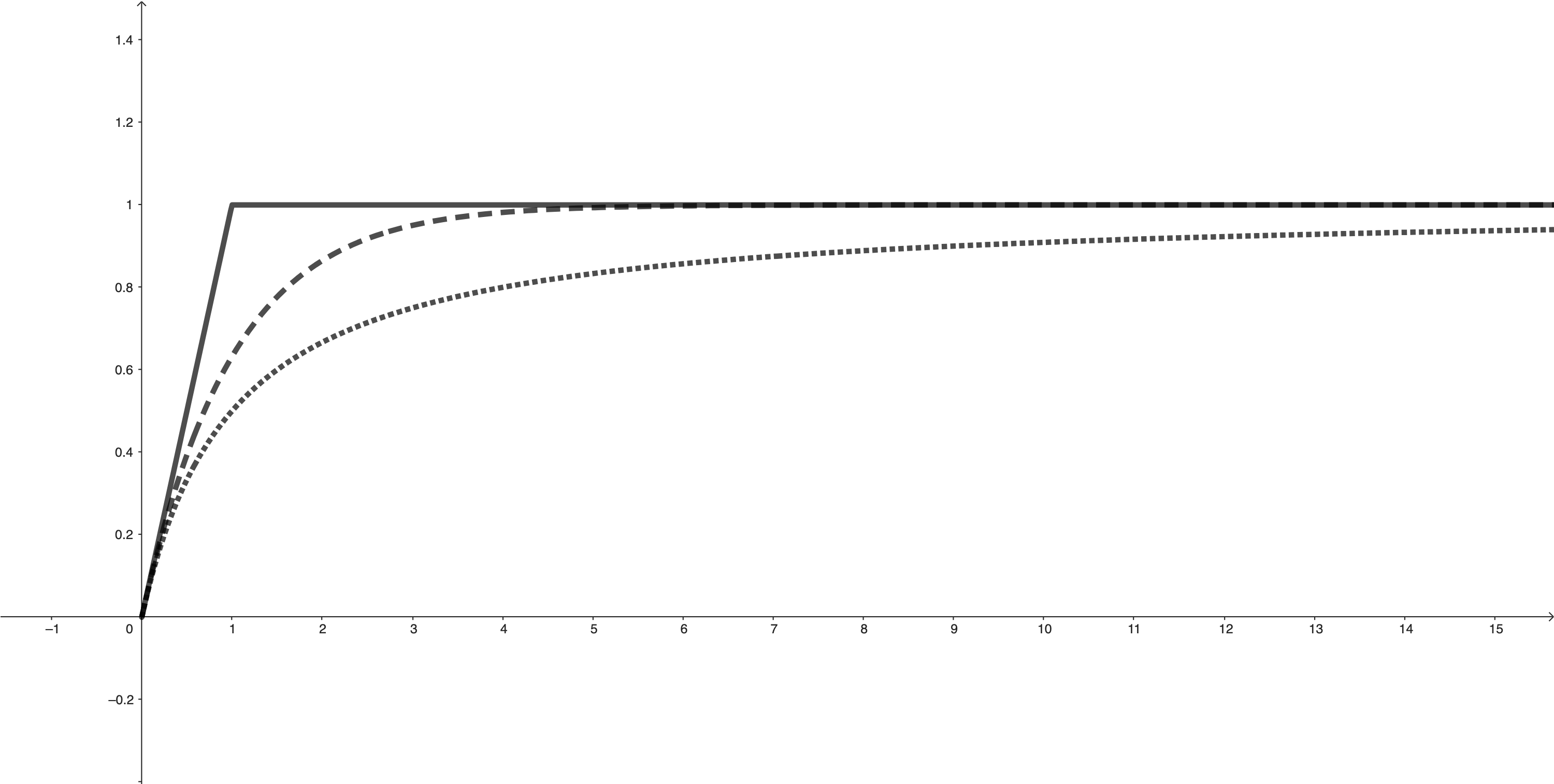}
    \caption{solid line: $\theta(t) = \min(1, t)$, dashed line: $\theta(t) = 1 - e^{-t}$, dotted line: $\theta(t) = \frac{t}{1+t}$}
    \label{Fig - Exemples}
\end{figure}

\begin{Example}
    A large number of maps $\theta$ satisfy the conditions of Corollary \ref{Coro - Propriété (P)}, including the functions used in \cite{COP-1}. For instance,
    \begin{itemize}[label=\sbt]
        \item $\theta(t) = \min(1, t)$, since $\theta(t) = t$ if $t < 1$ and $1-\theta(t) = 0$ if $t > 1$.
        \item $\theta(t) = \frac{t^\alpha}{1+t^\alpha}$, since $\theta(t) \le t^\alpha$ if $t < 1$ and $1-\theta(t) = \frac{1}{1+t^\alpha} \le t^{-\alpha}$ if $t > 1$.
        \item $\theta(t) = 1 - e^{-t}$, since $\frac{\theta(t)}{\sqrt{t}} \to 0$ if $t \to 0$, and $t(1-\theta(t)) = te^{-t} \to 0$ if $t \to +\infty$.
    \end{itemize}
    Graphs of these functions are displayed in Figure \ref{Fig - Exemples}.
\end{Example}

\section{Application : characterisation of linear isometries on some spaces}\label{Section 2}

In this section, we consider two examples of Fréchet spaces, for which we will describe the linear isometries.

\subsection{The space $\Hol(\D)$}

We consider here the space of holomorphic functions on the unit disc $\D$. 
We endow it with the family of supremum seminorms, defined by
\[ \norm{f}_{\infty, n} := \sup_{z \in K_n} \abs{f(z)}, \]
with $K_n = \{z \in \C : \abs{z} \le 1 - \frac{1}{n}\}$, giving an increasing and exhaustive sequence of compact subsets of $\D$.
In this case, the main theorem of \cite{COP-1} gives a complete characterisation of the linear isometries on the space $(\Hol(\D), \{\norm{\cdot}_{\infty, n}\}_{n \ge 1})$, hence for all the distances $d_{(\theta,r)}$, with $\theta$ satisfying the hypotheses of Corollary \ref{Coro - Propriété (P)}. Recall that $\T$ is the unit circle.

\begin{Theorem}
    Let $T : \Hol(\D) \to \Hol(\D)$ be a linear and continuous operator such that for all $f \in \Hol(\D)$ and $n \ge 1$,
    \[ \norm{T(f)}_{\infty, n} = \norm{f}_{\infty, n}. \]
    Then, there exist $\alpha, \beta \in \T$ such that for all $f \in \Hol(\D)$ and $z \in \D$,
    \[ T(f)(z) = \alpha f(\beta z) =: T_{\alpha, \beta}(f)(z). \]
\end{Theorem}

\begin{Remark}
    The proof of the theorem remains true if we consider only two seminorms. Hence, we may replace \guillemotleft{} $n \ge 1$ \guillemotright{} by \guillemotleft{} $n \in \{2, 3\}$ \guillemotright{}.
\end{Remark}

We may ask a natural question, and wonder what happens if we modify the family of seminorms, and consider the $H^p$-seminorms instead of the supremum ones. For $1 \le p < \infty$, let us define
\[ \norm{f}_{p, n} = \left(\frac{1}{2\pi} \int_0^{2\pi} \abs{f\left(\left(1 - \frac 1n\right) e^{i\theta}\right)}^p \dd \theta \right)^{1/p}. \]  
The following theorem gives an answer when the seminorms are not Hilbertian.

\begin{Theorem}
    Let $1 \le p < \infty$, $p \neq 2$ and $T : \Hol(\D) \to \Hol(\D)$ be a linear and continuous operator such that for all $f \in \Hol(\D)$ and $n \ge 1$,
    \[ \norm{T(f)}_{p, n} = \norm{f}_{p, n}. \]
    Then there exist $\alpha, \beta \in \T$ such that for all $f \in \Hol(\D)$ and $z \in \D$,
    \[ T(f)(z) = \alpha f(\beta z) =: T_{\alpha, \beta}(f)(z). \]
\end{Theorem}

In order to prove this result, we will need the next four results.

\begin{Proposition}[{\cite[Lemma 1]{Xiao-Zhu}}]\label{Prop 1}
    Let $0 < p \le \infty$, $0 \le r < 1$ and $f \in \Hol(\D)$ be non-constant. Set
    \[ M_{p,r}(f) = \left(\frac{1}{2\pi} \int_0^{2\pi} \abs{f\left(r e^{i\theta}\right)}^p \dd \theta \right)^{1/p}. \]
    Then, $r \mapsto M_{p,r}(f)$ is strictly increasing on $[0, 1)$.
\end{Proposition}

\begin{Proposition}[Forelli's theorem, \cite{Forelli}]\label{Prop 2}
    Let $1 \le p < \infty$, $p \neq 2$. Let $H^p(\D)$ the Hardy space of the unit disc, namely
    \[ H^p(\D) := \left\{f \in \Hol(\D) : \norm{f}_p = \lim_{r \to 1} \left(\frac{1}{2\pi} \int_0^{2\pi} \abs{f\left(r e^{i\theta}\right)}^p \dd \theta \right)^{1/p} < +\infty \right\}. \]
    It is a Banach space, endowed with the norm $\norm{\cdot}_p$.
    If an operator $T : H^p(\D) \to H^p(\D)$ is a linear isometry of $(H^p(\D), \norm \cdot _p)$, then there exist $m \in H^p(\D)$ and $\varphi : \D \to \D$ an inner function (that is $\abs{\varphi(re^{i\theta})} \to 1$ as $r \to 1$ for almost every $\theta \in \R$) such that
    \[ T(f) = m(f \circ \varphi), \quad f \in H^p(\D). \]
\end{Proposition}

\begin{Proposition}[\cite{Milner}]\label{Prop 3}
    Let $0 < r_1 < r_3$, and $f$ be an analytic map on  $\{z \in \C : r_1 \le \abs{z} \le r_3\}$. For $r \in [r_1, r_3]$, let $M(r) = \sup\{\abs{f(z)} : \abs{z} = r\}$. Then, $\log(M(r))$ is a convex function of $\log(r)$, i.e. for all $r_1 < r_2 < r_3$,
    \begin{equation*}
        \log\left(\frac{r_3}{r_1}\right) \log(M(r_2))
        \le \log\left(\frac{r_3}{r_2}\right) \log(M(r_1))
            + \log\left(\frac{r_2}{r_1}\right) \log(M(r_3)).
    \end{equation*}
    Moreover, the inequality is strict, unless $f(z) = c z^n$ for some $c \in \C$ and $n \in \Z$.
\end{Proposition}

\begin{Corollary}[{\cite[Corollary 3.2]{COP-2}}]\label{cor - Carac rot}
    Let $0 < r_1 < r_3$, and $f$ be an analytic map on the annulus $\{z \in \C : r_1 \le \abs{z} \le r_3\}$. Let $r_2 \in (r_1, r_3)$, and assume that $f(r_j \T) \subset r_j \T$ for all $j = 1, 2, 3$. Then, there exists $c \in \T$ such that $f(z) = c z$.
\end{Corollary}

\begin{proof}[Proof of Theorem 3]
    Let us divide the proof into three steps.
    
    \uline{Step 1}: We show that $T(\indic)$ is a constant and unimodular function. To do so, note that
    \[ \norm{T(\indic)}_{p,2} = \norm{\indic}_{p,2}
        = 1 = \norm{\indic}_{p,3} = \norm{T(\indic)}_{p,3}. \]
    If we write this equation with the definition of the seminorms, then
    \[ \left(\frac{1}{2\pi} \int_0^{2\pi} \abs{T\indic\left(\frac 12 e^{i\theta}\right)}^p \dd \theta \right)^{1/p}
        = \left(\frac{1}{2\pi} \int_0^{2\pi} \abs{T\indic\left(\frac 23 e^{i\theta}\right)}^p \dd \theta \right)^{1/p} = 1. \]
    Hence, $M_{p,1/2}(T\indic) = M_{p,2/3}(T\indic)$.
    This means that $r \mapsto M_{p,r}(T\indic)$ is not strictly increasing. Thus, by Proposition \ref{Prop 1}, $T(\indic)$ is a constant map. Denote by $\alpha \in \C$ the constant. Then,
    \[ M_{p,1/2}(T\indic) = \abs{\alpha} = 1. \medskip \]

    \uline{Step 2}: Let $n \ge 1$, and set $r_n = (n-1)/n$. If $\C[z]$ denotes the space of complex-valued polynomials, then $\C[z]$ is a dense subset of $H^p(r_n\D)$. Hence, for all $f \in H^p(r_n\D)$, there exists a sequence $(Q_k)_{k \ge 0} \subset \C[z]$ such that
    \[ \norm{f - Q_k}_{p,n} \xrightarrow[k \to \infty]{} 0. \]
    We define $\tilde T_n(f) = \underset{k \to \infty}{\lim} T(Q_k)$. Then, $\tilde T_n : H^p(r_n\D) \to H^p(r_n\D)$ has the following properties.
    \begin{itemize}[label=$\star$]
        \item The operator $\tilde T_n$ is well-defined. Indeed, if $(Q_k)$ and $(R_k)$ are two sequences of polynomials such that 
        \[ \lim_{k \to \infty} \norm{f - Q_k}_{p,n} 
            = \lim_{k \to \infty} \norm{f - R_k}_{p,n} = 0, \]
        then, we obtain
        \begin{align*}
            \norm{T(Q_k) - T(R_k)}_{p,n} = \norm{T(Q_k - R_k)}_{p,n}
            & = \norm{Q_k - R_k}_{p,n} \\ & \le \norm{Q_k - f}_{p,n} + \norm{f - R_k}_{p,n}
            \xrightarrow[k \to \infty]{} 0.
        \end{align*}

        \item Let $f \in H^p(r_n\D)$. Since $T$ is an isometry for $\norm \cdot _{p,n}$, we have
        \[ \norm{\tilde T_n(f)}_{p,n} = \lim_{k \to \infty} \norm{T(Q_k)}_{p,n}
        = \lim_{k \to \infty} \norm{Q_k}_{p,n} = \norm{f}_{p,n}. \]
        Thus, $\tilde T_n$ is an isometry of $H^p(r_n\D)$.
    \end{itemize}
    
    Using Forelli's theorem (Proposition \ref{Prop 2}), for all $n \ge 1$, there exist $m_n \in \Hol(r_n\D)$ and an inner function $\psi_n \in \Hol(r_n\D)$ such that for all $f \in H^p(r_n\D)$, if $\varphi_n(z) = r_n \psi_n(z/r_n)$, then
    \[ \tilde T_n(f) = m_n (f \circ \varphi_n). \medskip \]

    \uline{Step 3}: Note that if $Q \in \C[z]$, using the constant sequence $(Q_k = Q)_{k \ge 0}$, for all $n \ge 1$, 
    \[ T(Q) = \tilde T_n(Q) = m_n (Q \circ \varphi_n). \]
    We study two special polynomials.
    \begin{itemize}[label=$\star$]
        \item If $Q = \indic$, then $\tilde T_n(\indic) = m_n = T(\indic) \equiv \alpha \in \T$. \smallskip
        \item If $Q(z) = z$, then $\tilde T_n(z) = m_n \varphi_n = \alpha \varphi_n = T(z) \implies \varphi_n = \overline \alpha T(z)$. 
    \end{itemize}
    Finally, for all $Q \in \C[z]$, we can write $T(Q) = \alpha (Q \circ \varphi)$, with $\varphi(r_n \T) \subset r_n \T$ for all $n \ge 1$ (since $\varphi_n$ was inner).
    Using Corollary \ref{cor - Carac rot}, there exists $\beta \in \T$ such that $\varphi(z) = \beta z$. Hence,
    \[ T(Q)(z) = \alpha Q(\beta z) = T_{\alpha, \beta}(Q). \]

    To conclude, we only need to see that $\C[z]$ is a dense subset of $(\Hol(\D), \{\norm{\cdot}_{\infty, n}\}_{n \ge 0})$. Thus, for all $f \in \Hol(\D)$, $T(f) = T_{\alpha, \beta}(f)$.
\end{proof}

\begin{Remark}
    We may reduce the assumptions of the theorem, since only three different seminorms are sufficient to use Hadamard's three-circle theorem.
\end{Remark}

\begin{Remark}
    In \cite{Nov-Ob}, the authors obtained a characterisation of the linear isometries of the space $S^p = \{f \in \Hol(\D) : f' \in H^p(\D)\}$, endowed with the norm
    \[ \norm{f}_{S^p} := \sup_{z \in \D} \abs{f(z)} + \lim_{r \to 1} \left(\frac{1}{2\pi} \int_0^{2\pi} \abs{f\left(r e^{i\theta}\right)}^p \dd \theta \right)^{1/p}. \]
    They proved that, if $p \neq 2$, these are of the form
    \[ T_{\alpha, \beta} : f \mapsto \alpha f(\beta z), \]
    with $\alpha, \beta \in \C$ such that $\abs{\alpha} = \abs{\beta} = 1$. Therefore, with only one seminorm, one can show that the linear isometries of the space $\Hol(\D)$ endowed with the family of seminorms defined by
    \[ \norm{f}_{S^p, n} := \sup_{\abs z \le 1 - \frac 1n} \abs{f(z)} + \left(\frac{1}{2\pi} \int_0^{2\pi} \abs{f\left(\left(1 - \frac 1n\right) e^{i\theta}\right)}^p \dd \theta \right)^{1/p}. \]
    are the $T_{\alpha, \beta}$.
\end{Remark}

\subsection{The space $\mathcal C(U)$, surjective case}

We now focus on the space of all continuous complex-valued functions on an open set $U \subset \K^N$ ($\K = \R$ or $\C$), endowed with the supremum seminorms. We recall their definition:
\[ \norm{f}_{\infty, n} = \sup_{z \in K_n} \abs{f(z)}. \]
Here, we consider $(K_n)_{n \ge 0}$ a general increasing and exhaustive sequence of compact sets, i.e.
\[ K_n \subset K_{n+1} \quad \text{ and } \quad
    \bigcup_{n \ge 0} K_n = U. \]

Our starting point is Banach--Stone theorem (\cite[p.170]{Banach}, \cite[p.25]{Fleming-Jamison}), which characterises the linear isometries of $\mathcal C(K)$ onto $\mathcal C(Q)$ (endowed with the supremum norm), where $K$ and $Q$ are compact metric spaces. They first proved the real-valued case, but the complex-valued one follows by linearity.

\begin{Theorem}\label{Thm - Banach-Stone}
    Let $T : \mathcal C(K) \to \mathcal C(Q)$ be a linear, continuous and \textbf{surjective} operator such that for all $f \in \mathcal C(K)$,
    \[ \norm{T(f)}_{\infty, Q} := \sup_{z \in Q} \abs{T(f)(z)}
        = \sup_{z \in K} \abs{f(z)} =: \norm{f}_{\infty, K}. \]
    Then there exist a continuous unimodular function $h$ on $K$, and a homeomorphism $\varphi$ from $K$ onto $Q$ such that for all $f \in \mathcal C(K)$ and $z \in K$,
    \[ T(f)(z) = h(z) f(\varphi(z)). \]
\end{Theorem}

Thanks to this result, we can obtain a complete description of the linear isometries of the space $(\mathcal C(U), \{\norm{\cdot}_{\infty, n}\}_{n \ge 0})$.

\begin{Theorem}\label{Thm - Iso surj C(U)}
    Let $T : \mathcal C(U) \to \mathcal C(U)$ be a linear, continuous and \textbf{surjective} operator such that for all $f \in \mathcal C(U)$ and $n \ge 0$,
    \[ \norm{T(f)}_{\infty, n} = \norm{f}_{\infty, n}. \]
    Then there exist a continuous unimodular function $h$ on $U$ and a homeomorphism $\varphi$ from $U$ onto $U$ satisfying $\varphi(K_n) = K_n$ for all $n \ge 0$, such that for all $f \in \mathcal C(U)$ and $z \in U$,
    \[ T(f)(z) = h(z) f(\varphi(z)). \]
\end{Theorem}

\begin{proof}
    Let $n \ge 0$, and $T_n : \mathcal C(K_n) \to \mathcal C(K_n)$ be defined by
    \[ T_n(f) = (T \tilde f)_{|K_n}, \]
    where $\tilde f$ is a continuous extension of $f$ on $U$. Then, we have the following properties.
    \begin{itemize}
        \item The operator $T_n$ is well-defined. Indeed, if $f \in \mathcal C(K_n)$, and if $\tilde f$ and $\check{f}$ are two continuous extensions of $f$ on $U$, then
        \[ \norm{T\tilde f - T \check f}_{\infty, n}
            = \norm{\tilde f - \check f}_{\infty, n}
            = \norm{f - f}_{\infty, n} = 0, \]
        since $T$ is a $\norm{\cdot}_{\infty, n}$-isometry, and $\tilde f _{|K_n} = \check f _{|K_n} = f$. Therefore, $(T\tilde f)_{|K_n} = (T \check f)_{|K_n}$. \medskip

        \item The operator $T_n$ is a $\norm{\cdot}_{\infty, n}$-isometry. Indeed, for all $f \in \mathcal C(K_n)$,
        \[ \norm{T_n(f)}_{\infty, n} = \norm{T(\tilde f)}_{\infty, n}
            =  \norm{\tilde f}_{\infty, n} = \norm{f}_{\infty, n}. \]

        \item Finally, for all $n \ge 0$, $f \in \mathcal C(U)$ and $z \in K_n$, we have $T(f)(z) = T_n(f_{|K_n})(z)$.
    \end{itemize}

    Using Theorem \ref{Thm - Banach-Stone}, for all $n \ge 0$, there exist a continuous unimodular map $h_n$ on $K_n$, and a homeomorphism $\varphi_n$ from $K_n$ onto $K_n$ such that for all $f \in \mathcal C(U)$ et $z \in K_n$,
    \[ T(f)(z) = T_n(f_{|K_n})(z) = h_n(z) f(\varphi_n(z)). \]

    Let $n \ge 0$, and $z \in K_n$. We consider two particular maps.
    \begin{itemize}
        \item If $f = \indic$, we obtain $h_n(z) = T_n(\indic_{|K_n})(z) = T(\indic)(z)$. \smallskip
        \item If $f(z) = e_1(z) = z$, we obtain $\varphi_n(z) = \overline{h_n(z)} T(e_1)(z) = \overline{T(\indic)(z)} T(e_1)(z)$.
    \end{itemize}
    Hence, for all $z \in U$,
    \[ T(f)(z) = h(z) f(\varphi(z)), \]
    where $h = T(\indic)$ and $\varphi = \overline{T(\indic)} T(e_1)$. Moreover, $h$ is continuous on $U$, and for all $z \in U$, there exists $n_0 \ge 0$ such that $z \in K_{n_0}$. Thus, $\abs{h(z)} = \abs{h_{n_0}(z)} = 1$. In the same way, we show that
    \begin{itemize}
        \item $\varphi$ is continuous on $U$. Indeed, if $z \in U$, then there exists $n_0 \ge 0$ such that $z$ lies in the interior of $K_{n_0}$. Since $\varphi = \varphi_n$ on $K_{n_0}$ and $\varphi_n$ is continuous at $z$, $\varphi$ is continuous at $z$. \smallskip
        
        \item $\varphi$ is injective. Indeed, let $z_1, z_2 \in U$, and $n_0 \ge 0$ such that $z_1, z_2 \in K_{n_0}$. If $\varphi(z_1) = \varphi(z_2)$, then $\varphi_{n_0}(z_1) = \varphi_{n_0}(z_2)$, so $z_1 = z_2$ since $\varphi_{n_0}$ is injective. \smallskip

        \item $\varphi$ is surjective. Indeed, if $z_0 \in U$, there exists $n_0 \ge 0$ such that $z_0 \in K_{n_0}$. Since $\varphi_n$ is surjective, there exists $w_0 \in K_{n_0} \subset U$ such that $\varphi(w_0) = \varphi_{n_0}(w_0) = z_0$. \smallskip

        \item Using the same argument as for $\varphi$, $\varphi^{-1}$ is continuous on $U$.
    \end{itemize}
    Finally, $\varphi$ is a homeomorphism from $U$ onto $U$, and because $\varphi = \varphi_n$ on $K_n$ and $\varphi_n(K_n) = K_n$, we have $\varphi(K_n) = K_n$ for all $n \ge 0$.
\end{proof}

\begin{Example}
    Let us focus on two special open sets.
    \begin{enumerate}
        \item Let $U = (0, 1)$, and $K_n = [a_n, b_n]$, with $0 < a_0 \le b_0 < 1$, $(a_n)$ decreasing, $(b_n)$ increasing, $a_n \to 0$ and $b_n \to 1$. The main question is to describe the homeomorphisms $\varphi$ from $(0, 1)$ onto $(0, 1)$ such that $\varphi(K_n) = K_n$.
        Note that if $\varphi : (0, 1) \to (0, 1)$ is continuous and bijective, then it is strictly monotonous.
        \begin{itemize}
            \item If $\varphi$ is strictly increasing, then for all $n \ge 0$, $\varphi(a_n) = a_n$ and $\varphi(b_n) = b_n$.
            \item If $\varphi$ is strictly decreasing, then for all $n \ge 0$, $\varphi(a_n) = b_n$ and $\varphi(b_n) = a_n$.
        \end{itemize}
        Some examples of maps $\varphi$ are shown in Figure \ref{Fig - C(0,1)}.
    \end{enumerate}

    \begin{figure}[H]
        \centering
        \includegraphics[width=0.3\linewidth]{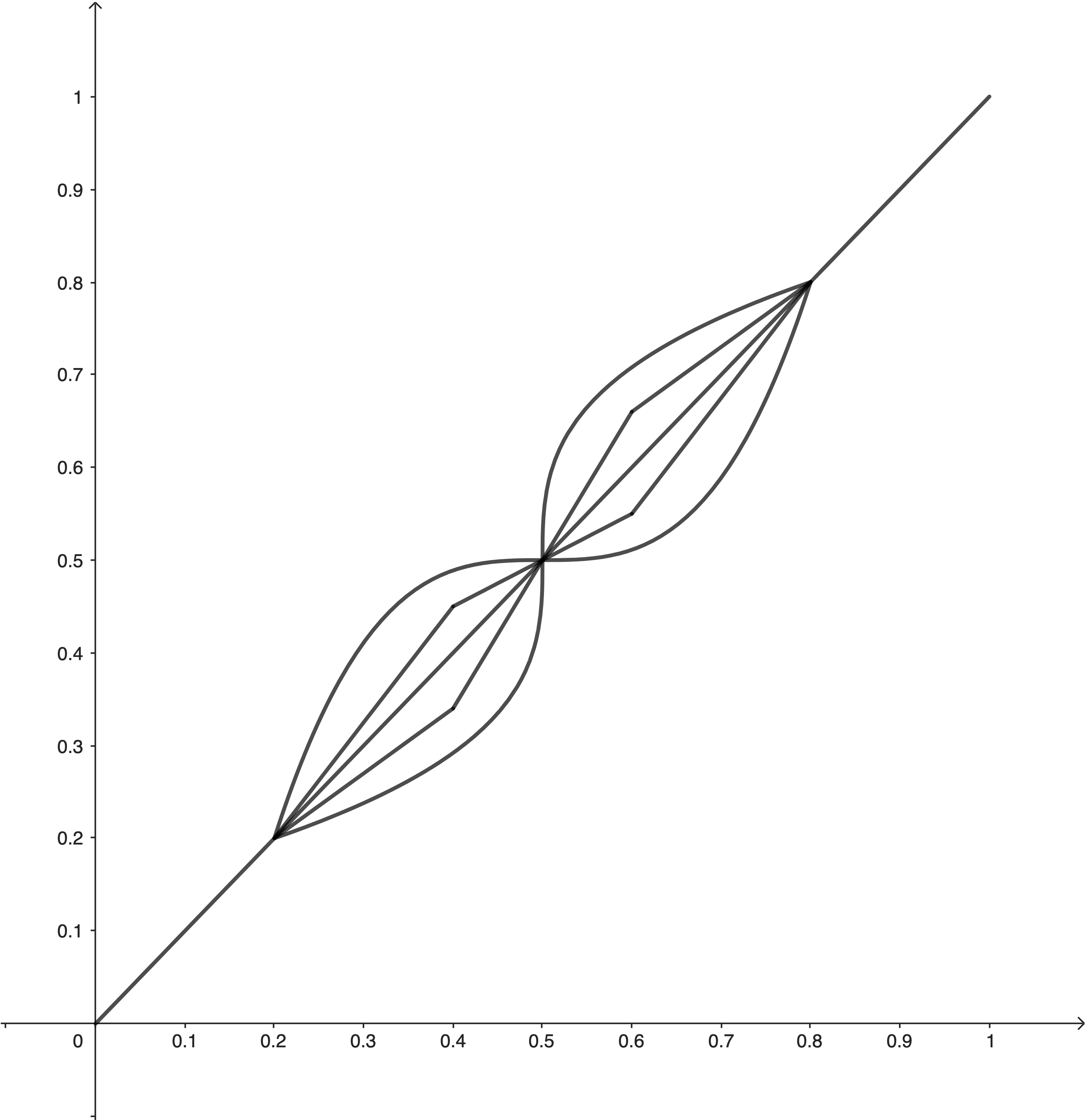}
        \hspace{1cm}
        \includegraphics[width=0.3\linewidth]{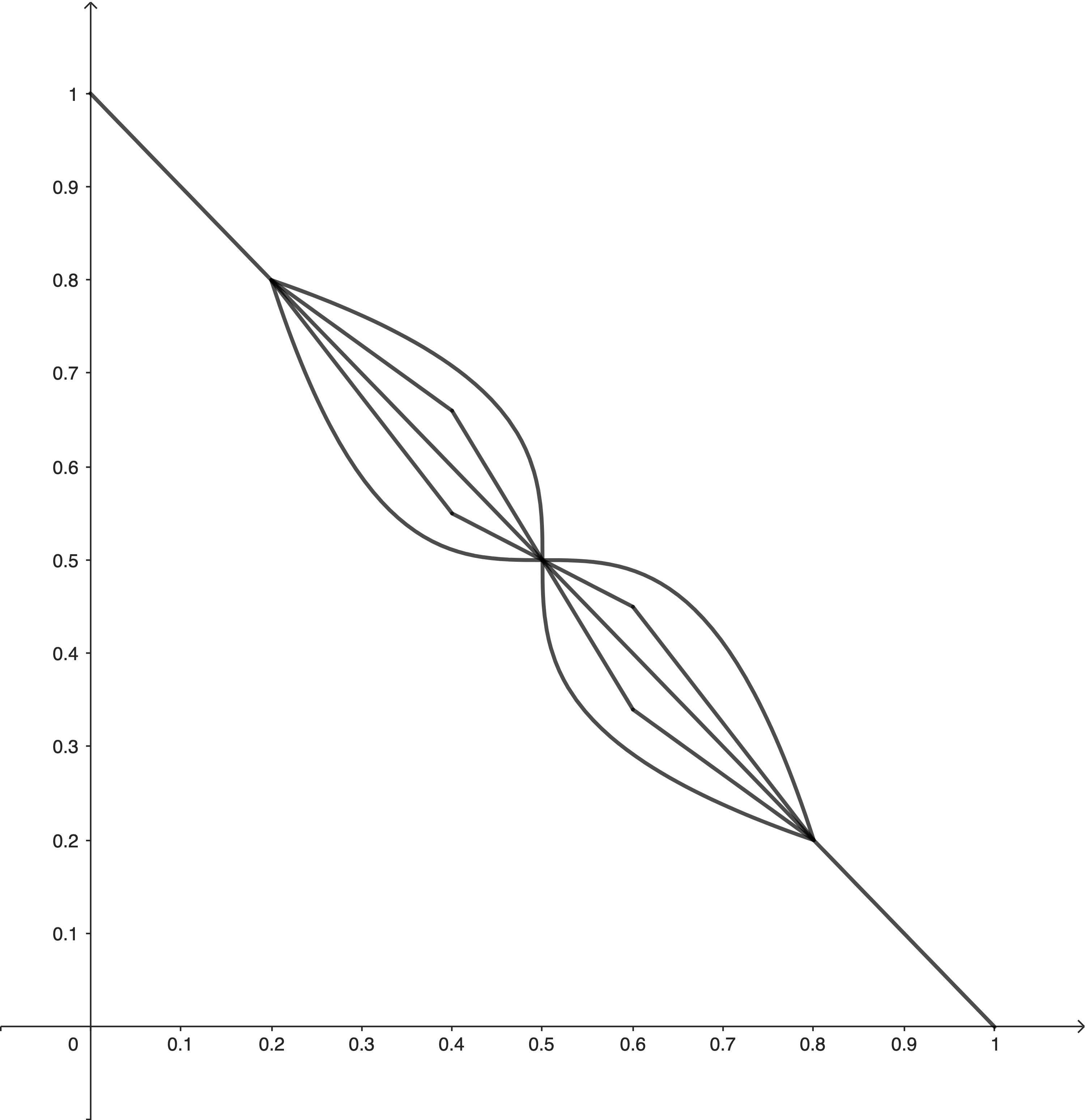}
        \caption{Homeomorphisms $\varphi : (0, 1) \to (0, 1)$ s.t. $\varphi(\frac 12) = \frac 12$ and $\varphi([\frac{1}{5}, \frac{4}{5}]) = [\frac{1}{5}, \frac{4}{5}]$}
        \label{Fig - C(0,1)}
    \end{figure}

    \begin{enumerate}
        \item[2.] Let $U = \D$, and $K_n = r_n \D$, with $0 \le r_0 < 1$, $(r_n)$ increasing and $r_n \to 1$. It is way more difficult to obtain a characterisation of the homeomorphisms $\varphi$ from $\D$ onto $\D$ such that $\varphi(K_n) = K_n$.
        We know that these $\varphi$ must satisfy, for all $n \ge 0$,
        \[ \varphi(r_n \T) = r_n \T \quad \text{ and } \quad
            \varphi(\{r_n \le \abs{z} \le r_{n+1}\}) = \{r_n \le \abs{z} \le r_{n+1}\}. \]
        However, inside the annuli $\{r_n \le \abs{z} \le r_{n+1}\}$, the behaviour is free (as long as it is continuous). In Figure \ref{Fig - C(D)}, we can see one example of $\varphi$.
    \end{enumerate}

    \begin{figure}[H]
        \centering
        \includegraphics[width=0.3\linewidth]{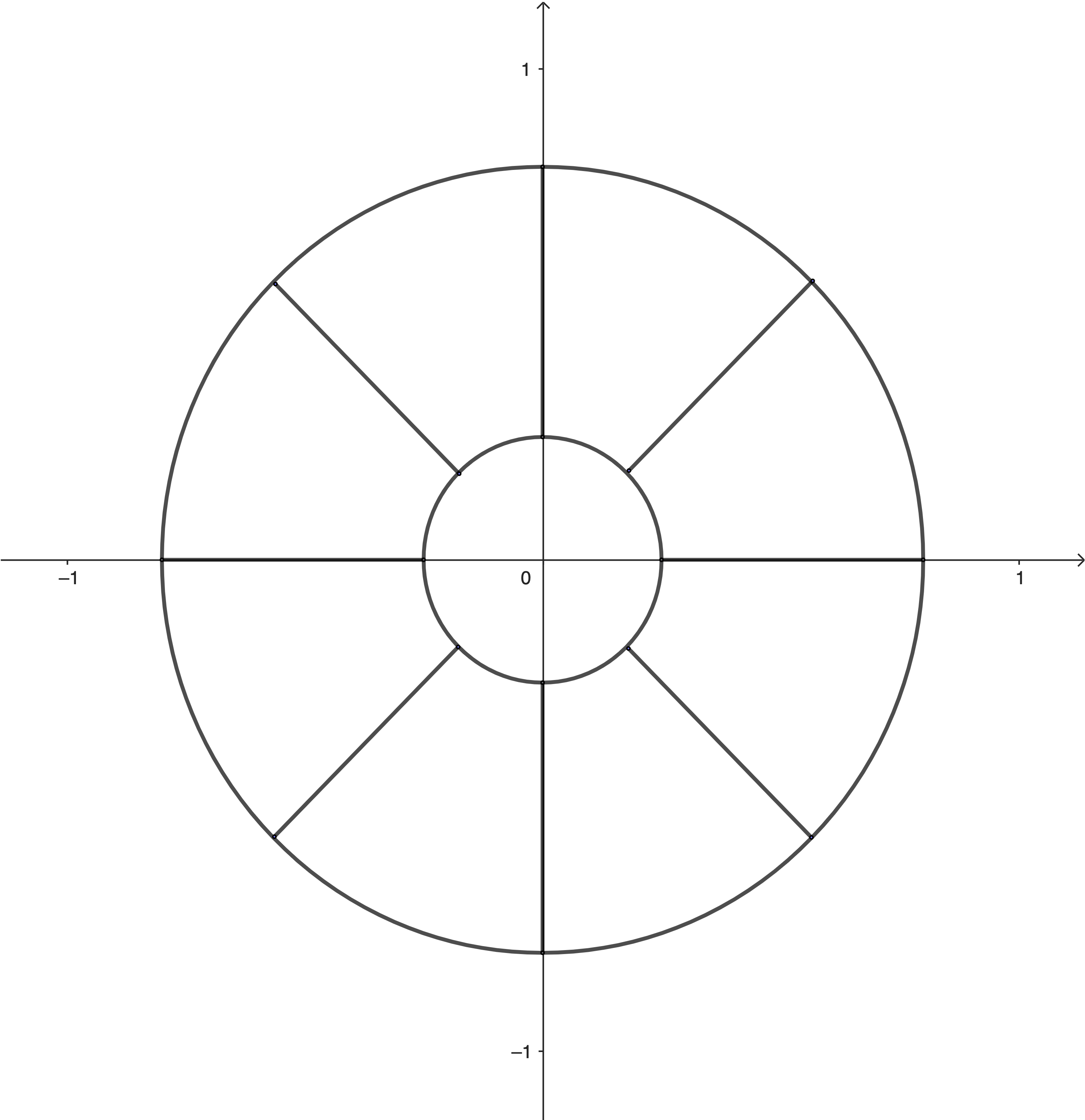}
        \begin{minipage}[b]{0.1\textwidth}
        \centering $\xrightarrow[]{\; \varphi \;}$ \\ {\color{white} . \\ . \\ . \\ . \\ .} \end{minipage}
        \includegraphics[width=0.3\linewidth]{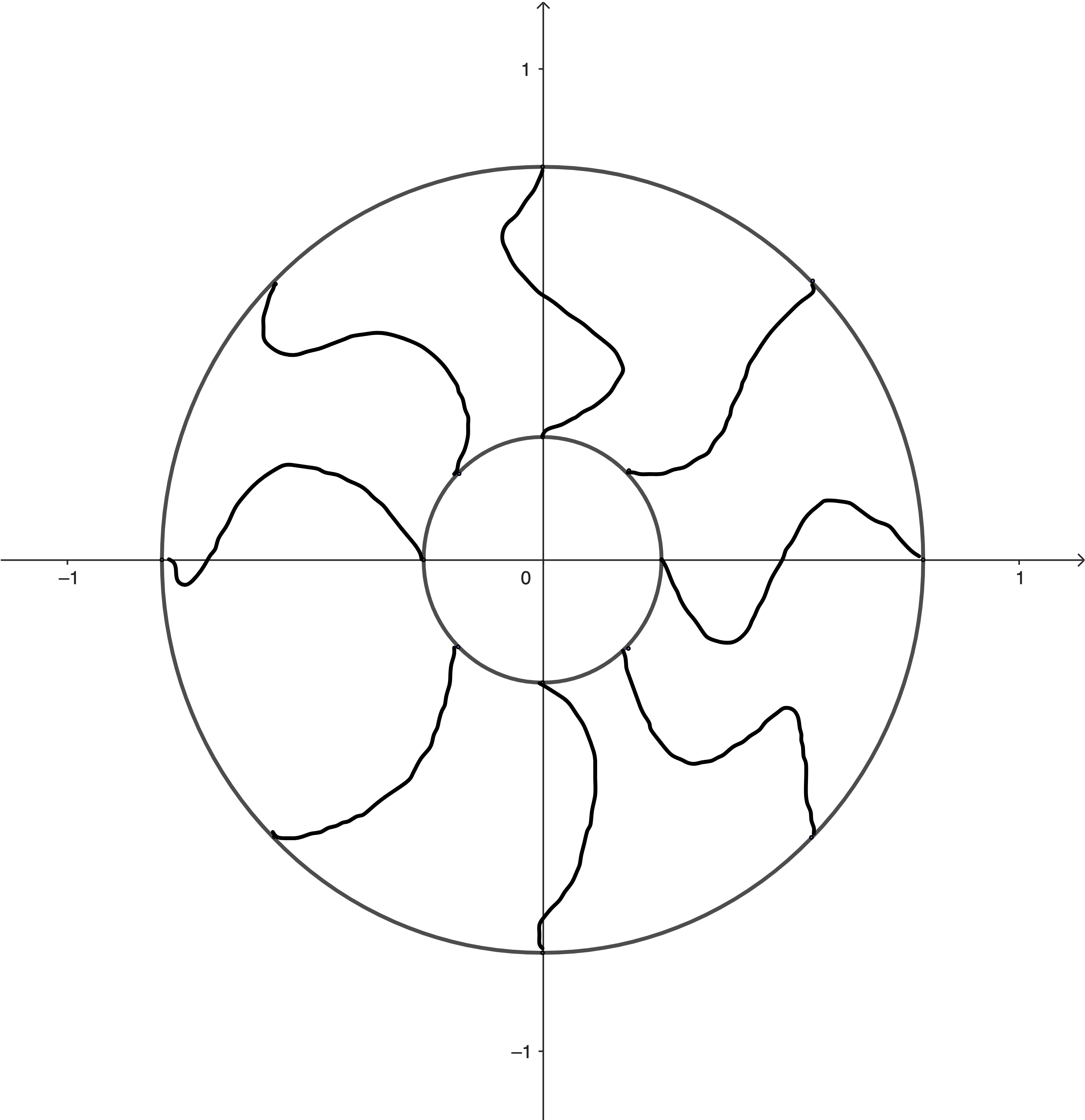}
        \caption{Homeomorphism $\varphi : \D \to \D$ s.t. $\varphi(K_n) = K_n$, with $K_0 = \frac 14\D$ and $K_1 = \frac 45\D$ (we assume that $\varphi$ is the identity map outside the annulus)}
        \label{Fig - C(D)}
    \end{figure}
\end{Example}

\subsection{The space $\mathcal C(U)$, nonsurjective case}

When we remove the surjectivity assumption, we need to use a generalisation of the Banach--Stone theorem, proved by Holszty\'nski \cite{Hol}
(see also \cite{Nov} and \cite[Theorem 2.3.10]{Fleming-Jamison}).

\begin{Theorem}\label{Thm - Novinger}
    Let $K$ be a compact set, and $T : \mathcal C(K) \to \mathcal C(K)$ a linear operator that is a $\norm{\cdot}_{\infty, K}$-isometry. Then there exists a closed subset $L \subset K$, a continuous map $h : L \to \T$, and a continuous surjective map $\varphi : L \to K$, such that for all $f \in \mathcal C (K)$ and $z \in L$,
    \[ T(f)(z) = h(z) f(\varphi(z)). \]
    Moreover, if we set 
    \[ h = T(\indic_K) \quad \text{ and } \quad
       \Phi(x)(f) = \overline{h(x)}T(f)(x), \]
    it is possible to extend $h$ continuously  on $K$ (keeping the unimodular property), and express $\varphi$ through a continuous map $\Phi : K \to \mathcal C(K)^*$, such that if $z \in L$ and $f \in \mathcal C(K)$, we have $\Phi(z)(f) = f(\varphi(z))$. Finally, for all $z \in K$, $\norm{\Phi(z)}^* = 1$, and
    \[ T(f)(z) = h(z) [\Phi(z)(f)], 
        \quad f \in \mathcal C(K), \quad z \in K. \]
\end{Theorem}

This result will be the key point to obtain a description of the nonsurjective linear isometries of $\mathcal C(U)$.

\begin{Theorem}\label{Thm - Iso C(U)}
    Let $T : \mathcal C(U) \to \mathcal C(U)$ be a linear continuous operator such that for all $f \in \mathcal C(U)$ and $n \ge 0$, 
    \[ \norm{T(f)}_{\infty, n} = \norm{f}_{\infty, n}. \]
    Then there exist two continuous maps $h : U \to \T$ and $\Phi : U \to \mathcal C(U)^*$ such that for all $f \in \mathcal C(U)$ and $z \in U$, $\abs{\Phi(z)(f)} \le \min \{\norm{f}_{\infty, n} : n \ge 0, \; z \in K_n\}$, and
    \[ T(f)(z) = h(z) [\Phi(z)(f)]. \]
\end{Theorem}

\begin{proof}
    We start in the same way as in the proof of Theorem \ref{Thm - Iso surj C(U)}. For $n \ge 0$, we define the operator $T_n : \mathcal C(K_n) \to \mathcal C(K_n)$ by
    \[ T_n(f) = (T \tilde f)_{|K_n}, \]
    where $\tilde f$ is a continuous extension of $f$ on $U$. Then, $T_n$ is well-defined, is a $\norm{\cdot}_{\infty, n}$-isometry and for $z \in K_n$, $T(f)(z) = T_n(f_{|K_n})(z)$. Hence, by Theorem \ref{Thm - Novinger}, for all $n \ge 0$, $f \in \mathcal C(U)$ and $z \in K_n$, we can write
    \[ T(f)(z) = T_n(f_{|K_n})(z) = h_n(z) [\Phi_n(z)(f_{|K_n})], \]
    with $h_n : K_n \to \T$ and $\Phi_n : K_n \to \mathcal C(K_n)^*$ continuous on $K_n$. Moreover, $h_n$ is defined as $T_n(\indic_{K_n})$, so if $h = T(\indic)$, we have $h$ continuous on $U$, and $h = h_n$ on $K_n$. Since $h_n$ is unimodular on $K_n$, $h$ is unimodular on $U$. Therefore,
    \[ T(f)(z) = T_n(f_{|K_n})(z) = h(z) [\Phi_n(z)(f_{|K_n})]. \]

    Let $z \in U$, and set $\Phi(z)(f) = \Phi_n(z)(f_{|K_n})$ if $z \in K_n$. Then,
    \begin{itemize}
        \item $\Phi$ is well-defined, because if $m \ge n$ and $z \in K_n \subset K_m$, we have
        \begin{align*}
            \Phi_n(z)(f_{|K_n}) 
            = \overline{h_n(z)} T_n(f_{|K_n})(z)
            & = \overline{h(z)} T(f)(z) \\
            & = \overline{h_m(z)} T_m(f_{|K_m})(z)
            = \Phi_m(z)(f_{|K_m}).
        \end{align*} 

        \item $\Phi(z)$ is a linear form on $\mathcal C(U)$, because if $f,g \in \mathcal C(U)$ and $\lambda \in \C$, for ${n_0} \ge 0$ such that $z \in K_{n_0}$, we have $\Phi(z)(f) = \Phi_{n_0}(z)(f_{|K_{n_0}}) \in \C$, and
        \begin{align*}
            \Phi(z)(f+\lambda g)
            & = \Phi_{n_0}(z)((f+\lambda g)_{|K_{n_0}}) \\
            & = \Phi_{n_0}(z)(f_{|K_{n_0}}) + \lambda \Phi_n(g_{|K_{n_0}})
            = \Phi(z)(f) + \lambda \Phi(z)(g).
        \end{align*}

        \item $\Phi(z)$ is continuous. Indeed, if $f \in \mathcal C(U)$, for $n_0 \ge 0$ such that $z \in K_{n_0}$,
        \[ \abs{\Phi(z)(f)} = \abs{\Phi_{n_0}(z)(f_{|K_{n_0}})}
            \le \norm{\Phi_{n_0}(z)}^* \norm{f}_{\infty, n_0}
            = \norm{f}_{\infty, n_0}. \]
    \end{itemize}
    To conclude, for $z \in U$ and $f \in \mathcal C(U)$, we have
    \[ T(f)(z) = h(z) [\Phi(z)(f)]. \qedhere \]
\end{proof}

\section{Conclusion}

This work has provided a very general result on isometries on Fr\'echet spaces, and applied it to the study of various standard spaces.
It is expected that it will provide a motivation for future research, and some developments that have been suggested include the study of weighted spaces of differentiable functions, the stability of isometries under perturbations, and the corresponding theory of nonlinear mappings.
\smallskip

\noindent \textbf{Acknowledgments:}
The authors are very grateful to Robert Eymard for some suggestions which motivated this work. 
This research is partly supported by the Bézout Labex, funded by ANR, reference ANR-10-LABX-58.

\end{document}